\documentclass[11pt,a4paper]{article}

\usepackage{amsmath,amssymb}

\newcommand{\diam}{\operatorname{diam}}

\newcommand{\counte}{theorem}
\newtheorem{theorem} {\bf Theorem}[section]

\newtheorem{defn}[\counte]{\bf Definition}

\newtheorem{lemma}[\counte]{\bf Lemma}

\newtheorem{coro}[\counte]{\bf Corollary}

\newtheorem{remark}[\counte]{\bf Remark}

\allowdisplaybreaks

\renewcommand{\thefootnote}{\fnsymbol{footnote}}

\begin{document}

\renewcommand{\thefootnote}{\arabic{footnote}}

\centerline{\bf ON ALMOST ISOMETRY THEOREM IN ALEXANDROV SPACES }
\centerline{\bf WITH CURVATURE BOUNDED BELOW \footnote{Supported  by
NSFC 10801011 and 10826052. \hfill{$\,$}}$^{, }$
\footnote{Mathematics Subject Classification (2000):
53-C20.\hfill{$\,$}}}

\vskip3mm

\centerline{ Xiaole Su$^{*}$, Hongwei Sun$^{\dagger}$, Yusheng
Wang$^{\ddagger,}$ \footnote{The corresponding author (E-mail:
wwyusheng@gmail.com). \hfill{$\,$}}}

{\it \noindent$^{*,\ \ddagger}$ Mathematics Department, Beijing
Normal University, Beijing, 100875 P.R.C.

\noindent$^{\dagger}$ Mathematics Department, Capital Normal
University, Beijing, 100037 P.R.C.}

\vskip2mm

\centerline{\bf Abstract.}

In this paper we give a new proof for an almost isometry theorem in
Alexandrov spaces with curvature bounded below.

\vskip1mm


\vskip6mm

Due to the great work by Perel'man on Poincar\'e conjecture,
Alexandrov geometry (especially with curvature bounded below)
together with Gromov-Hausdorff convergence theory attracts more and
more attentions.

A fundamental and significant work on Alexandrov spaces with
curvature bounded below is of Burago-Gromov-Perel'man (\cite{BGP}).
One important result in \cite{BGP} is an almost isometry theorem
(see Theorem 0.1 below). We find a key lemma of its proof is
incorrect (see ``Lemma'' 1.2 and Example 1.3 below). We suppose that
the authors of \cite{BGP} missed some condition. In the present
paper we adjust the conditions of the lemma so that the conclusion
of it still holds (see Lemma 2.1 below). Unfortunately, from the
modified lemma the original proof of the theorem cannot go through.
For this reason, we supply a new proof for the theorem in this
paper.


\setcounter{section}{-1}

\section{Notations and main theorem}

We first give some notations, which are almost copied from
\cite{BGP}.

$\bullet$ $|xy|$ always denotes the distance between two points $x$
and $y$ in a metric space.

$\bullet$ For any three points $p,q,r$ in a length space, we
associate a triangle $\triangle\tilde p\tilde q\tilde r$ on the
$k$-plane (2-dimensional complete and simply-connected Riemannian
manifold of constant curvature $k$) with $|\tilde p\tilde
q|=|pq|,|\tilde p\tilde r|=|pr|$ and $|\tilde r\tilde q|=|rq|$. For
$k\leqslant 0$ and for $k>0$ with $|pq|+|pr|+|qr|\leqslant
2\pi/\sqrt k$, such a triangle always exists. We denote by $\tilde
\angle pqr $ the angle of the triangle $\triangle\tilde p\tilde
q\tilde r$ at vertex $\tilde q$.

$\bullet$ $M$ always denotes an Alexandrov space with curvature
bounded below by $k$, which is a length space and in which there
exists a neighborhood $U_x$ around any $x\in M$ such that for any
four (distinct) points $(a; b,c,d)$ in $U_x$
$$\tilde \angle bac+\tilde \angle bad+\tilde \angle cad\leqslant
2\pi.$$

$\bullet$ A point $p\in M$ is called an $(n, \delta)$-strained point
if there are $n$ pairs of points $(a_i,b_i)$ distinct from $p$ such
that for $i\ne j$
$$\begin{aligned} \tilde \angle a_ipb_i>\pi-\delta,\  \tilde \angle
a_ipa_j>\pi/2-\delta,& \\
 \tilde \angle a_ipb_j>\pi/2-\delta,\ \tilde \angle
 b_ipb_j>\pi/2-\delta.
\end{aligned} $$
$\{(a_i,b_i)\}_{i=1}^n$ is called an $(n, \delta )$-strainer at $p$
(which is obviously a generalization of a coordinate frame). We say
that the $(n, \delta )$-strainer $\{(a_i,b_i)\}_{i=1}^n$ at $p$ is
$R$-long if $|a_ip|>\dfrac{R}{\delta}$ and
$|b_ip|>\dfrac{R}{\delta}$ for all $i$. And we denote by $M(n,
\delta, R)$ the set of points with $R$-long $(n,\delta)$-strainer in
$M$.

$\bullet$ An important fact is that if any neighborhood of a point
$p\in M$ contains an $(n,\delta)$-strained point ($\delta$ is
sufficient small) but no $(n+1,\delta)$-strained point, then any
neighborhood of any other point in $M$ has the same property (see \S
6 in \cite{BGP}). And it follows that the dimension of such $M$ is
defined to be $n$.

$\bullet$ We always denote by $\varkappa(\cdot)$ or
$\varkappa(\cdot,\cdot)$ (resp. $C$) a positive function which is
infinitesimal at zero (e.g.
$\varkappa(\delta,\delta_1)\longrightarrow0$ as $\delta,
\delta_1\longrightarrow0$) (resp. a constant depending only on $n$);
however we do not distinguish any two distinct $\varkappa$-functions
with the same parameters (resp. any two such constants) when we use
them.

$\bullet$ A map $f$ between metric spaces $(X, d_1)$ and $(Y, d_2)$
is called a GH$_\epsilon$-approximation if $B_\epsilon(f(X))=Y$ and
$|d_2(f(x_1),f(x_2))-d_1(x_1, x_2)|<\epsilon$ for any $x_1, x_2\in
X$.

$\bullet$ $f:(X, d_1)\longrightarrow(Y, d_2)$ is called a {\it
$\varkappa(\delta)$-almost distance-preserving map} if
$$\left|1-\dfrac{|f(x)f(y)|}{|xy|}\right|
<\varkappa(\delta) \text { for any } x, y\in X;$$ and if in addition
$f$ is a bijection, $f$ is called a {\it $\varkappa(\delta)$-almost
isometry}.

$\bullet$ We say that $\bar f:(X, d_1)\longrightarrow(Y, d_2)$ is
$\nu$-close to $f$ if $|f(x)\bar f(x)|<\nu$ for any $x\in X$.

\vskip2mm

Now we formulate the almost isometry theorem in \cite{BGP} mentioned
at the beginning.

{\noindent \bf Theorem 0.1 (Theorem 9.8 in \cite{BGP})} { \it Let
$M_1$ and $M_2$ be two compact $n$-dimensional Alexandrov spaces
with the same low curvature bound, and let $h: M_1 \to M_2$ be a
GH$_\nu$-approximation. Then for sufficiently small $\delta$ and
$\dfrac{\nu}{R\delta^3}$, there exists a
$\varkappa(\delta,\frac{\nu}{R\delta^3})$-almost distance preserving
map $\overline{h}: M_1(n,\delta, R)\to M_2$ which is $C\nu$-close to
$h$.}

\vskip1mm

It is not difficult to conclude from Theorem 0.1 the following
important corollary.

\vskip1mm

{\noindent \bf Corollary 0.2 (\cite{BGP})} {\it In Theorem 0.1, if
in addition each point of $M_2$ is $(n,\delta)$-strained, then there
exists a $\varkappa(\delta, \nu)$-almost isometry $\overline{h}:
M_1\longrightarrow M_2$ which is $C\nu$-close to $h$.}

\vskip1mm

Theorem 0.1 (or Corollary 0.2) plays an important role when one
studies a converging sequence (with respect to the Gromov-Hausdorff
distance defined by the
GH-appro-\\
ximation) of $n$-dimensional Alexandrov spaces with the same low
curvature bound.

In this paper we give the proof of the following sharper version of
Theorem 0.1.

{\noindent \bf Theorem A} { \it Let $M_1$ and $M_2$ be two compact
$n$-dimensional Alexandrov spaces with the same low curvature bound,
and let $h: M_1 \longrightarrow M_2$ be a GH$_\nu$-approximation.
Then for sufficiently small $\delta$ and $\nu<\delta^2 R$, there
exists a $\varkappa(\delta)$-almost distance preserving map
$\overline{h}: M_1(n,\delta, R)\longrightarrow M_2$ which is
$C\nu$-close to $h$.}

\vskip2mm

The construction of $\bar h$ is almost copied from \cite{BGP} (see
Section 3). The main difference between our proof and the proof of
Theorem 0.1 in \cite{BGP} is how to verify that $\bar h$ almost
preserves distance (see Section 4). Of course, we use some ideas and
results in \cite{BGP}.

\vskip2mm

\noindent{\bf Remark 0.3 }\  In \cite{Ya} Yamaguchi  proved that,
without the assumption of the dimension of $M_1$, there is an almost
Lipschitz submersion from $M_1$ to $M_2$ if each point of $M_2$ is
$(n,\delta)$-strained in Theorem 0.1. This result (which appears as
a conjecture in \cite{BGP}) coincides with Corollary 0.2 if the
dimension of $M_1$ is $n$. The key approach to construct the almost
Lipschitz submersion in \cite{Ya} is to embed an Alexandrov space
with curvature bound below $M$ into $L^2(M)$. Compared with it, the
base of the construction of $\bar h$ of Theorem 0.1 (or A) is that
$M_1(n,\delta, R)$ is locally almost isometric to the
$n$-dimensional Euclidean space (see Theorem 1.1 below).
\endremark


\section{Center of mass and a key lemma in \cite{BGP}}

The main tool in the construction of $\bar h$ (\cite{BGP}) in
Theorem 0.1 (or A) is ``center of mass''. Recall that the center of
mass of a set of points $Q=\{q_1,q_2,\cdots, q_l\}\subset
\mathbb{R}^n$ with weights $W=(w_1,w_2, \cdots, w_l)$ (where
$\sum_{j=1}^l w_j=1$ and $w_j\geqslant 0$) is defined to be
$$Q_W=\sum\limits_{j=1}^l w_j q_j.$$

The construction of the center of mass for a set of points in $M$ is
based on the following important result.

\begin{theorem}{\bf(Theorem 9.4 in \cite{BGP})}
Let $M$ be an $n$-dimensional Alexandrov space with curvature
bounded below, and let $\{(a_i,b_i)\}_{i=1}^n$ be an
$(n,\delta)$-strainer at $p\in M$. Then there exist neighborhoods
$U$ and $V$ around $p$ and $(|a_1p|,|a_2p|,\cdots,|a_np|)\in
\mathbb{R}^n$ respectively such that $$f: U\longrightarrow
V\subset\mathbb{R}^n \text{ given by }
f(q)=(|a_1q|,|a_2q|,\cdots,|a_nq|)$$ is a $\varkappa(\delta,
\delta_1)$-almost isometry, where $\delta_1=\max\limits_{1\leqslant
i\leqslant n}\{|pa_i|^{-1}, |pb_i|^{-1}\}\cdot\diam U.$
\end{theorem}

If $Q=\{q_1,q_2,\cdots, q_l\}$ belongs to $U$ in Theorem 1.1, and if
in addition $f(U)$ is convex in $\mathbb{R}^n$, then the center of
mass of $Q$ with weights $W$ is defined to be (\cite{BGP})
$$Q_W=f^{-1}\left(\sum\limits_{j=1}^l w_j f(q_j)\right).$$
Obviously, $Q_W$ depends on the choice of the $(n,\delta)$-strainer
at $p$.

Now we give the key lemma in \cite{BGP} mentioned at the beginning
of Section 0, which plays a crucial role in verifying that $\bar h$
in Theorem 0.1 almost preserves distance.

\vskip1mm

\noindent{\bf ``Lemma'' 1.2}  {\it Let $p, U, \{(a_i,b_i)\}_{i=1}^n$
and $f$ be the same as in Theorem 1.1, and let
$\{(s_i,t_i)\}_{i=1}^n$ be another $(n,\delta)$-strainer at $p$ with
$$
\delta_1=\max\left\{\dfrac{\diam U}{\min_i\{|pa_i|,|pb_i|\}},
\dfrac{\max_i\{|pa_i|,|pb_i|\}}{\min_i\{|ps_i|,|pt_i|\}}\right\}.
$$
And let $Q=\{q_1,\cdots, q_l\}$ and $R=\{r_1,\cdots,r_l\}$ be two
sets of points in $U$ with
$$
\max_j\{|q_jr_j|\}< (1+\delta)\min_j\{|q_jr_j|\}, \hbox{ and }
|\max_j\tilde \angle s_iq_jr_j-\min_j\tilde\angle s_iq_jr_j|< \delta
\hbox{ for } i=1,\cdots n.
$$
Assume that $f(U)$ is convex in $\mathbb{R}^n$. Then for any weights
$W^1$ and $W^2$ such that $||W^1-W^2||<\delta_1$, the centers of
mass $Q_{W^1}$ and $R_{W^2}$ (with respect to the strainer
$\{(a_i,b_i)\}$) satisfy that
$$
\left|1-\dfrac{|q_jr_j|}{|Q_{W^1}
R_{W^2}|}\right|<\varkappa(\delta,\delta_1)$$ and $|\tilde \angle
s_iq_jr_j-\tilde\angle s_iQ_{W^1}R_{W^2}|<
\varkappa(\delta,\delta_1)\hbox{ for } j=1,\cdots,l \hbox{ and }
i=1,\cdots,n.$}

\vskip1mm

Due to the following counterexample, we don't think that this lemma
is correct.

\vskip1mm

\noindent{\bf Example 1.3 } 1. In fact, if $q_j=r_j$ for
$j=1,\cdots,l$ and $W^1\ne W^2$, then $Q_{W^1}\ne R_{W^2}$ and thus
$$\left|1-\dfrac{|q_jr_j|}{|Q_{W^1}
R_{W^2}|}\right|=|1-0|>\varkappa(\delta,\delta_1).$$

\noindent2. If $|q_jr_j|\ll||W^1-W^2||$ for all $j$, then
``$\left|1-\dfrac{|q_jr_j|}{|Q_{W^1}
R_{W^2}|}\right|<\varkappa(\delta,\delta_1)$'' does not hold.

\vskip2mm

Inspired by the example, we add some stronger restriction on the
weights $W^1$ and $W^2$ (see Lemma 2.1 below) so that the conclusion
in ``Lemma 1.2'' still holds.


\section{Modified key lemma}

In this section we give a modified version of ``Lemma'' 1.2 which is
formulated as follows (for convenience we divide it into two parts).

\begin{lemma} \label{lemma2.6}
Let $p, U, \{(a_i,b_i)\}_{i=1}^n$ and $f$ be the same as in Theorem
1.1, and let $\{(s_i,t_i)\}_{i=1}^n$ be another
$(n,\delta)$-strainer at $p$ with
$$
\max\left\{\dfrac{\diam U}{\min_i\{|pa_i|,|pb_i|\}},
\dfrac{\max_i\{|pa_i|,|pb_i|\}}{\min_i\{|ps_i|,|pt_i|\}}\right\}<\delta.
$$
And let $Q=\{q_1,\cdots, q_l\}$ and $R=\{r_1,\cdots,r_l\}$ be two
sets of points in $U$. Then the following conclusions hold.

\noindent(2.1.1) The following statements are equivalent:

$(1)$ $ | \tilde \angle a_iq_jr_j-\tilde\angle a_iq_{j'}r_{j'}|<
\varkappa(\delta) $ for $i=1,\cdots,n$;

$(2)$ $| \tilde \angle s_iq_jr_j- \tilde\angle s_iq_{j'}r_{j'}|<
\varkappa(\delta) $ for $i=1,\cdots,n$.

\noindent(2.1.2) Assume that $f(U)$ is convex in $\mathbb{R}^n$, and
assume that
\begin{equation}\label{eq2.2}
\max_j\{|q_jr_j|\}< (1+\varkappa(\delta))\min_j\{|q_jr_j|\} \text{
and }\end{equation}
\begin{equation}\label{eq2.3}
|\max_j\tilde \angle a_iq_jr_j-\min_j\tilde\angle
a_iq_jr_j|<\varkappa(\delta) \text{ for } i=1,\cdots,n.
\end{equation}
Then for any weights $W^1$ and $W^2$ with $||W^1-W^2||\cdot
\max\limits_{j,j'}|r_jr_{j'}|<
{\varkappa(\delta)\min\limits_j|q_jr_j|}$, the centers of mass
$Q_{W^1}$ and $R_{W^2}$ (with respect to the strainer
$\{(a_i,b_i)\}$) satisfy that
$$
\left|1-\dfrac{|q_jr_j|}{|Q_{W^1} R_{W^2}|}\right|<\varkappa(\delta)
\hbox{ and }
$$
$$
|\tilde \angle a_iq_jr_j-\tilde\angle a_iQ_{W^1}R_{W^2}|<
\varkappa(\delta) \hbox{ for } j=1,\cdots,l \hbox{ and }
i=1,\cdots,n.
$$
\end{lemma}

(2.1.1) is proved in [1] (for convenience of readers we give its
proof in Appendix). For the proof of (2.1.2) we need Lemmas 2.2 and
2.4.
\begin{lemma} \label{lemma2.4}Let $Q=\{q_1,q_2,\cdots, q_l\}$ and $R=\{r_1,r_2,\cdots,
r_l\}$ be two sets of points in $\mathbb{R}^n$, and let
$W^i=(w_1^i,w_2^i, \cdots, w_l^i)$ be two weights with $i=1,2$. Then
$$
\overrightarrow{Q_{W^1}R_{W^2}}=\sum_{j=1}^l
w_j^1\overrightarrow{q_jr_j}+\sum_{j=1}^l(w_j^2-w_j^1)\overrightarrow{r_{j_0}r_j},
$$
for any $j_0\in \{1,2,\cdots,l\}$.
\end{lemma}

\vskip2mm

\noindent{\it Proof}.\ \ Straightforward computation gives
$$
\begin{aligned} \hskip3.5cm \overrightarrow{Q_{W^1}R_{W_2}} &
 = \sum\limits_{j=1}^l w^2_j r_j -\sum\limits_{j=1}^l w^1_j q_j\\
  &
 =\sum\limits_{j=1}^l w^1_j  (r_j-q_j)+\sum\limits_{j=1}^l (w^2_j-w^1_j) r_j \\
  &
 = \sum\limits_{j=1}^l w^1_j \overrightarrow{q_jr_j}+\sum\limits_{j=1}^l (w^2_j-w^1_j) r_j
 -\sum\limits_{j=1}^l (w^2_j-w^1_j) r_{j_0} \\
& =\sum_{j=1}^l
w_j^1\overrightarrow{q_jr_j}+\sum_{j=1}^l(w_j^2-w_j^1)\overrightarrow{r_{j_0}r_j}.
\hskip3.5cm\Box\end{aligned}$$

To simplify further considerations, we use the following definition.

\begin{defn}{\rm
For sets of points $Q=\{q_1,q_2\}$ and $R=\{r_1,r_2\}$ in
$\mathbb{R}^n$, we say that $\overrightarrow{q_1r_1}$ is {\it
$\varkappa(\delta)$-almost parallel} to $\overrightarrow{q_2r_2}$ if
$$\angle (\overrightarrow{q_1r_1},
\overrightarrow{q_2r_2})<\varkappa(\delta);$$ and if in addition
$$
\left|1-\dfrac{|q_1r_1|}{|q_2r_2|}\right|<\varkappa(\delta),
$$
we say that $\overrightarrow{q_1r_1}$ is {\it
$\varkappa(\delta)$-almost equal} to $\overrightarrow{q_2r_2}$.}
\end{defn}

\begin{lemma} \label{lemma2.7} Let $p, U, \{(a_i,b_i)\}_{i=1}^n$ and $f$ be the same as
in Lemma \ref{lemma2.6}. Then for any points $x_1,x_2,y_1,y_2\in U$,
the following statements are equivalent:

\noindent$(1)$ $|\tilde \angle a_ix_1y_1-\tilde\angle a_ix_2y_{2}|<
\varkappa(\delta)$ for $i=1,2,\cdots,n$;

\noindent$(2)$ $\overrightarrow{f(x_1)f(y_1)}$ is
$\varkappa(\delta)$-almost parallel to
$\overrightarrow{f(x_2)f(y_2)}$.
\end{lemma}

Lemma 2.4 is implied in [1] (we will give its proof in Appendix).

\vspace{2mm}

\noindent{\bf Proof of (2.1.2).} According to Theorem 1.1 and Lemma
\ref{lemma2.7}, inequalities (\ref{eq2.2}) and (\ref{eq2.3}) imply
that $\overrightarrow{f(q_j)f(r_j)}$ are $\varkappa(\delta)$-almost
equal each other for $j=1,2,\cdots,l$. Therefore it follows from
Lemma \ref{lemma2.4} that $\overrightarrow{f(Q_{W^1})f(R_{W^2})}$ is
$\varkappa(\delta)$-almost equal to $\overrightarrow{f(q_j)f(r_j)}$
for every $j$ (note that $f(Q_{W^1})=\sum\limits_{j=1}^l w^1_j
f(q_j)$ and $f(R_{W^2})=\sum\limits_{j=1}^l w^2_j f(r_j)$, and
$||W^1-W^2||\cdot \max\limits_{j,j'}|r_jr_{j'}|<
{\varkappa(\delta)\min\limits_j|q_jr_j|}$). And thus the conclusion
in (2.1.2) follows from Lemma \ref{lemma2.7} and the fact that $f$
is a $\varkappa(\delta)$-almost isometry. \hfill$\Box$

\vskip2mm

At the end of this section we give a corollary of (2.1.1), which
will be used in gluing local almost isometries to a global one (see
next section).

\begin{coro} \label{coro2.10}
Let $p, U, \{(a_i,b_i)\}_{i=1}^n$ and $\{(s_i,t_i)\}_{i=1}^n$ be the
same as in Lemma \ref{lemma2.6}. Let $\{(a_i',b_i')\}_{i=1}^n$ be an
$(n,\delta)$-strainer at another point $p'$, and let $U'$ be a
neighborhood around $p$ determined by Theorem 1.1 (with respect to
$\{(a_i',b_i')\}$). Moreover we assume that $\{(s_i,t_i)\}_{i=1}^n$
is also an $(n,\delta)$-strainer at $p'$, and
$$
\left\{\dfrac{\diam U'}{\min_i\{|p'a'_i|,|p'b'_i|\}},\
\dfrac{\max_i\{|p'a'_i|,|p'b'_i|\}}{\min_i\{|p's_i|,|p't_i|\}}\right\}<\delta.
$$
Then for any points $x_1,x_2, y_1,y_2\in U_1\cap U_2$, the following
statements are equivalent:

\noindent$(1)$ $| \tilde \angle a_ix_1y_1- \tilde\angle
a_ix_2y_{2}|< \varkappa(\delta) $ for $i=1,\cdots,n$;

\noindent$(2)$ $| \tilde \angle a'_ix_1y_1- \tilde\angle
a'_ix_2y_{2}|< \varkappa(\delta) $ for $i=1,\cdots,n$.
\end{coro}


\setcounter{equation}{0}
\section{The construction of $\bar h$ in Theorem A}

In this section, we give the construction of the map $\bar h$ in
Theorem A, which is almost copied from [1].

\vskip2mm

Since the closure of $M_1(n,\delta, R)$ is compact, we can select
$x_j\in M_1(n,\delta, R)$ with $j=1,\cdots,N_1$ such that
\begin{equation}\label{eq3.*}\bigcup_{j=1}^{N_1}B_{x_j}(\delta R)\supset
\bigcup_{j=1}^{N_1}B_{x_j}(\frac13\delta R)\supset M_1(n,\delta,
R).\end{equation} Without loss of generality, we can assume that the
multiplicity of the cover $\{B_{x_j}(\delta R)\}$ is bounded by a
number $N$ depending only on the dimension $n$ (see Theorem 1.1 for
the dimension).

Since $x_j\in M_1(n,\delta, R)$, there exists an $R$-long
$(n,\delta)$-strainer $\{(s_i^j, t_i^j)\}_{i=1}^n$ at $x_j$ (with
$\min_i\{|x_js_i^j|,|x_jt_i^j|\}>\frac{R}{\delta}$), and thus there
exists a $\delta R$-long $(n,\delta)$-strainer $\{(a_i^j,
b_i^j)\}_{i=1}^n$ at $x_j$ (with
$\min_i\{|x_ja_i^j|,|x_jb_i^j|\}>R$) such that
$$
\dfrac{\max_i\{|x_ja_i^j|,|x_jb_i^j|\}}{\min_i\{|x_js_i^j|,|x_jt_i^j|\}}<\delta.
$$
Denote by $f_j$ and $U_j$ the associated map and the neighborhood
around $x_j$ in Theorem 1.1 with respect to the strainer $\{(a_i^j,
b_i^j)\}$. Moreover we select $U_j$ such that $f_j(U_j)$ is convex
in $\mathbb{R}^n$; and such that
\begin{equation}\label{eq3.**}B_{x_j}(2\delta R/3)\subset U_j\subset
B_{x_j}(\delta R)\end{equation} which implies  that
$$f_j\big|_{U_j} \text{ is a $\varkappa(\delta)$-almost isometry (see Theorem 1.1)}.$$

Since $h$ is a GH$_\nu$-approximation with $\nu< R\delta^2$,
$\{(h(a_i^j), h(b_i^j))\}_{i=1}^n$ and $\{(h(s_i^j),
h(t_i^j))\}_{i=1}^n$ are $(n,2\delta)$-strainers at $h(x_j)$. We
consider the associated map $g_j$ around $h(x_j)$ in Theorem 1.1
with respect to the strainer $\{(h(a_i^j), h(b_i^j))\}$, and we have
that
$$g_j^{-1}\big|_{f_j(U_j)} \text{ is a $\varkappa(\delta)$-almost
isometry}.$$

Obviously, $$h_j=g_j^{-1}\circ f_j \text{ is a
$\varkappa(\delta)$-almost isometry on each } U_j,$$ and for any
$x\in U_j$
$$\begin{aligned}
&|h_j(x)h(x)|=(1+\varkappa(\delta))|g_j(h_j(x))g_j(h(x))|\\
=&(1+\varkappa(\delta))|f_j(x)g_j(h(x))|\\
=&(1+\varkappa(\delta))\sqrt{(|a_1^jx|-|h(a_1^j)h(x)|)^2+\cdots+
(|a_n^jx|-|h(a_n^j)h(x)|)^2}\\
<&(1+\varkappa(\delta))\sqrt{n}\nu \ \text{ (note that $h$ is a
GH$_\nu$-approximation)},
\end{aligned}$$
i.e. each $h_j$ is $C\nu$-close to $h$ on $U_j$.

We will use center of mass to glue all these local almost isometries
$h_j$ to a global one. We first define weight functions\footnote{The
original definition in [1] is $\phi_j(x)=(1-2|xx_j|/(\delta R))^N$
if $x\in B_{x_j}(\delta R/2)$, but we find that power $1$ is
sufficient. A basic reason for this is that we only need Lipschitz
condition.} $\phi_j: M_1\longrightarrow \mathbb{R}$ by
$$
\phi_j(x)=\begin{cases}
1-\frac{2|xx_j|}{\delta R}, & x\in B_{x_j}(\delta R/2),\\
0, & x\in M_1\backslash B_{x_j}(\delta R/2).
\end{cases}
$$
Then for an arbitrary point $z\in M_1(n,\delta, R)$ we define a
sequence $\{z_j\}_{j=1}^{N_1}\subset M_2$ :
$$z_{j}=\begin{cases}
g_{j}^{-1}\left(\dfrac{\Sigma_{j-1}(z)}{\Sigma_{j}(z)}g_{j}(z_{j-1})+
\dfrac{\phi_{j}(z)}{\Sigma_{j}(z)}g_{j}(h_{j}(z))\right) & z\in U_{j}\\
z_{j-1}, & z\not\in U_{j}
\end{cases}, $$
where $z_0=h(z)$, $\Sigma_0(z)=0$ and
$\Sigma_{j}(z)=\sum\limits_{l=1}^{j}\phi_l(z)$ for $j\geqslant1$. A
basic fact is that
\begin{equation}\label{eq3.2}
\Sigma_{N_1}(z)>\frac13\ \text{ (see (\ref{eq3.*}))}.
\end{equation}

Now we define the desired map $\bar h:
M_1(n,\delta,R)\longrightarrow M_2$ in Theorem A by
$$\bar h(z)=z_{N_1} \hbox{ for any } z\in M_1(n,\delta,R).$$
Since each $h_j$ is $C\nu$-close to $h$, it is easy to see that
\begin{equation}\label{eq3.3}
|h_j(z)h_{j'}(z)|< C\nu \text{ and } |z_jh_{j'}(z)|< C\nu,
\end{equation}
and thus $$\bar h \text{ is } C\nu\text{-close to } h.$$

\vskip2mm

In next section we will {\bf verify that $\bar h$
$\varkappa(\delta)$-almost preserves distance}.


\setcounter{equation}{0}

\section{Verifying that $\bar h$ almost preserves distance}

In this section, we verify that $\bar h$ constructed in Section 3
almost preserves distance, i.e. for any $y, z\in M_1(n,\delta, R)$,
\begin{equation}\label{eq4.1}
\left|1-\dfrac{|\overline{h}(y)\overline{h}(z)|}{|yz|}\right|<\varkappa(\delta)
\text{ or }
\left||\overline{h}(y)\overline{h}(z)|-|yz|\right|<\varkappa(\delta)|yz|,
\end{equation}
and thus {\bf the proof of Theorem A is completed}.

We first observe that we only need to consider the case
``$|yz|<R\delta^{3/2}$''. In fact, if $|yz|\geqslant R\delta^{3/2}$,
then $||\overline{h}(y)\overline{h}(z)|-|yz||< C\nu<
CR\delta^2<|yz|\varkappa(\delta)$ (i.e. (\ref{eq4.1}) holds) because
$\overline{h}$ is $C\nu$-close to $h$ which is a
GH$_{\nu}$-approximation.

Without loss of generality, we assume that
$\phi_j(y)+\phi_j(z)\neq0$ for $1\leqslant j\leqslant N_2$, but
$\phi_j(y)+\phi_j(z)=0$ for $N_2<j\leqslant N_1$. Note that if
$\phi_j(y)\neq0$ (i.e., $y\in B_{x_j}(\delta R/2)$), then $z\in
B_{x_j}(\delta 2R/3)\subset U_j$ (see (\ref{eq3.**})) because $|yz|<
R\delta^{3/2}$ ($\delta$ is sufficient small). Then
$$y, z\in U_j \text{ for } j=1,\cdots, N_2 \text{ and }
y, z\not\in U_j \text{ for } j>N_2,$$ which implies that
$N_2\leqslant N$ (a number depending on $n$) and that
\begin{equation}\label{eq4.2}
\bar{h}(y)=y_{N_2} \text{ and } \bar{h}(z)=z_{N_2}.
\end{equation}
And we can define two new sequences $\{\overline{y}_j\}_{j=1}^{N_2}$
and $\{\overline{z}_j\}_{j=1}^{N_2}$ in $M_2$ (which are not
introduced in [1]):
$$\overline{y}_{j}=g_{j}^{-1}\left(\sum\limits_{l=1}^{j}
\dfrac{\phi_{l}(y)}{\Sigma_{j}(y)}g_{j}(h_l(y))\right) \text{ and }
\overline{z}_{j}=g_{j}^{-1}\left(\sum\limits_{l=1}^{j}
\dfrac{\phi_{l}(z)}{\Sigma_{j}(z)}g_{j}(h_l(z))\right).$$ Note that
\begin{equation}\label{eq4.3}
\overline{y}_j=y_j \text{ and } \overline{z}_j=z_j\text{ for }
j=1,2.
\end{equation}

Now we give two claims.

\noindent{\bf Claim 1}\footnote{The present proof is mainly inspired
by this observation.}:
$$
\left||\overline{y}_{N_2}\overline{z}_{N_2}|-|yz|\right|
<\varkappa(\delta)|yz|.
$$

\noindent{\bf Claim 2}:
$$
||y_{N_2}z_{N_2}|-|\overline{y}_{N_2}\overline{z}_{N_2}||<
\varkappa(\delta)|yz|.
$$

Obviously, {\bf Claims 1} and {\bf 2} (together with (\ref{eq4.2}))
imply (\ref{eq4.1}). Hence {\bf we only need to verify Claims 1 and
2}.

\vskip2mm

\noindent$\bullet$ {\bf The proof of  Claim 1}:

\vskip2mm

Note that $\overline{y}_{N_2}$ (resp. $\overline{z}_{N_2}$) is the
center of mass of $\{h_j(y)\}_{j=1}^{N_2}$ (resp.
$\{h_j(z)\}_{j=1}^{N_2}$) with weights
$W_y=(\frac{\phi_{1}(y)}{\Sigma_{N_2}(y)},\cdots,\frac{\phi_{N_2}(y)}{\Sigma_{N_2}(y)})$
(resp.
$W_z=(\frac{\phi_{1}(z)}{\Sigma_{N_2}(z)},\cdots,\frac{\phi_{N_2}(z)}{\Sigma_{N_2}(z)})$)
with respect to the $(n,\delta)$-strainer $\{(h(a_i^{N_2}),
h(b_i^{N_2}))\}_{i=1}^n$ at $h(x_{N_2})$. Then according to (2.1.2),
{\bf Claim 1} follows from the following three properties.

(i) Since each $h_j$ is a $\varkappa(\delta)$-almost isometry, we
have
\begin{equation}\label{eq4.?}
\max_{j}\{|h_j(y)h_{j}(z)|\}<(1+\varkappa(\delta))\min_{j}\{|h_j(y)h_{j}(z)|\}.
\end{equation}

(ii) For any fixed $j$,
\begin{equation}\label{eq4.6}
|\max_{l}\tilde \angle h(a_i^j)h_{l}(y)h_{l}(z)-\min_{l}\tilde\angle
h(a_i^j)h_{l}(y)h_{l}(z)|<\varkappa(\delta) \text{ for }
i=1,\cdots,n.
\end{equation}
This is proved in [1] (we give its proof in Appendix in which the
strainers $\{(s_i^j, t_i^j)\}$ will be used).

(iii) \begin{equation}\label{eq4.7}||W_y-W_z||\cdot
\max\limits_{j,j'}|h_j(z)h_{j'}(z)|<
{\varkappa(\delta)\min\limits_j|h_j(y)h_j(z)|}.
\end{equation}
In order to prove inequality (\ref{eq4.7}), we first give an
estimate
\begin{equation}\label{eq4.8}
\left|\dfrac{\phi_l(y)}{\Sigma_j(y)}-\dfrac{\phi_l(z)}{\Sigma_j(z)}\right|
\leqslant \dfrac{C|yz|}{\delta R\Sigma_j(y)}  \text{ for }
1\leqslant l\leqslant j\leqslant N_2.
\end{equation}
In fact, for any $1\leqslant l\leqslant N_2$ we have
$|\phi_l(y)-\phi_l(z)|=2\dfrac{||zx_l|-|yx_l||}{\delta R}\leqslant
\dfrac{2|yz|}{\delta R}$, and thus
$$
\begin{aligned}
\left|\dfrac{\phi_l(y)}{\Sigma_j(y)}-\dfrac{\phi_l(z)}{\Sigma_j(z)}\right|=&
\dfrac{1}{\Sigma_j(y)}\left|\phi_l(y)-\dfrac{\phi_l(z)
\Sigma_j(y)}{\Sigma_j(z)}\right|\\
= & \dfrac{1}{\Sigma_j(y)}\left|\phi_l(y)-\phi_l(z)-\phi_l(z)\dfrac{
\Sigma_j(y)-\Sigma_j(z)}{\Sigma_j(z)}\right|\\
\ \leqslant &
\dfrac{1}{\Sigma_j(y)}\max\limits_l\{\left|\phi_l(y)-\phi_l(z)\right|\}\cdot
N_2\\
\leqslant  & \dfrac{C|yz|}{\delta R\Sigma_j(y)}.
\end{aligned}
$$
Note that inequality (\ref{eq4.7}) follows from (\ref{eq4.8}),
$\Sigma_{N_2}(y)>\frac{1}{3}$ (see (\ref{eq3.2})) and
$|h_{j}(z)h_{j'}(z)|<C\nu< CR\delta^2$ (see (\ref{eq3.3})).

\vskip2mm

\noindent$\bullet$ {\bf The proof of  Claim 2}:

\vskip2mm

Put $\overrightarrow{\alpha}_j
=\overrightarrow{g_j({y}_{j})g_j(z_{j})}-
\overrightarrow{g_j(\overline{y}_{j})g_j(\overline{z}_{j})}$,
$j=1,\cdots,N_2$. Since each $g_j$ is a $\varkappa(\delta)$-almost
isometry, {\bf Claim 2} is equivalent to
\begin{equation}\label{eq4.9}
|\overrightarrow{\alpha}_{N_2}|< \varkappa(\delta)|yz|.
\end{equation}

\noindent{\bf Subclaim}:
\begin{equation}\label{eq4.10}
|\overrightarrow{\alpha}_j|\leqslant \dfrac{C|yz|\nu}{\delta
R\Sigma_j(y)} +\dfrac{\Sigma_{j-1}(y)}{\Sigma_j(y)}
(1+\varkappa(\delta))
|\overrightarrow{\alpha}_{j-1}|+\varkappa(\delta)|yz| \text{ for }
j=2,\cdots,N_2.
\end{equation}
It follows from the subclaim that
$$
\begin{aligned} |\overrightarrow{\alpha}_{N_2}| &  \leqslant
\dfrac{C|yz|\nu}{\delta R\Sigma_{N_2}(y)}+\varkappa(\delta)|yz|
+\dfrac{\Sigma_{N_2-1}(y)}{\Sigma_{N_2}(y)} (1+\varkappa(\delta))
|\overrightarrow{\alpha}_{N_2-1}|\\
& \leqslant \dfrac{C|yz|\nu}{\delta R\Sigma_{N_2}(y)}
+\varkappa(\delta)|yz|+\dfrac{\Sigma_{{N_2}-2}(y)}{\Sigma_{N_2}(y)}
(1+\varkappa(\delta))|\overrightarrow{\alpha}_{{N_2}-2}|\\
& \leqslant\ \ \cdots \\
& \leqslant \dfrac{C|yz|\nu}{\delta R\Sigma_{N_2}(y)}
+\varkappa(\delta)|yz|+\dfrac{\Sigma_{2}(y)}{\Sigma_{N_2}(y)}
(1+\varkappa(\delta))|\overrightarrow{\alpha}_{2}|\\
& < \varkappa(\delta)|yz|\ \ \text{(note that
$\Sigma_{N_2}(y)>\frac13,\ \nu< R\delta^2$ and
$|\overrightarrow{\alpha}_{2}|=0$ (see (\ref{eq4.3}))).}
\end{aligned}
$$

{\bf Now we only need to verify the subclaim}.

To simplify notations in the following computations, we let $\tilde
x$ denote  $g_j(x)$ for any $x\in U_j$.

Recall that
$$\tilde y_j=\dfrac{\Sigma_{j-1}(y)}{\Sigma_{j}(y)}\tilde y_{j-1}+
\dfrac{\phi_j(y)}{\Sigma_{j}(y)}\widetilde{h_j(y)}\text{\ \ and\ \ }
\widetilde{\overline{y}}_{j}=\sum_{l=1}^j\dfrac{\phi_l(y)}{\Sigma_j(y)}\widetilde{h_l(y)}$$
($\tilde z_j$ and $\widetilde{\overline{z}}_{j}$ have the same form
respectively). Through straightforward computation, one can get
$$
\begin{aligned}
 \overrightarrow{\alpha}_j
=&
\dfrac{\Sigma_{j-1}(y)}{\Sigma_{j}(y)}\left(\overrightarrow{\widetilde
y_{j-1}\widetilde z_{j-1}}
-\sum_{l=1}^{j-1}\dfrac{\phi_l(y)}{\Sigma_{j-1}(y)}\overrightarrow{\widetilde{h_l(y)}\widetilde{h_l(z)}}\right)
+\sum_{l=1}^{j-1}
\left(\dfrac{\phi_l(z)}{\Sigma_j(z)}-\dfrac{\phi_l(y)}{\Sigma_j(y)}\right)
\overrightarrow{\widetilde{h_l(z)}\widetilde{z}_{j-1}}
\end{aligned}
$$
Put $\overrightarrow{\beta}=\overrightarrow{\widetilde
y_{j-1}\widetilde
z_{j-1}}-\sum\limits_{l=1}^{j-1}\dfrac{\phi_l(y)}{\Sigma_{j-1}(y)}
\overrightarrow{\widetilde{h_l(y)}\widetilde{h_l(z)}}$ and
$\overrightarrow{\gamma}=\sum\limits_{l=1}^{j-1}
\left(\dfrac{\phi_l(z)}{\Sigma_j(z)}-\dfrac{\phi_l(y)}{\Sigma_j(y)}\right)
\overrightarrow{\widetilde{h_l(z)}\widetilde{z}_{j-1}},$ and thus
\begin{equation}\label{eq4.11}
\overrightarrow{\alpha}_j=\dfrac{\Sigma_{j-1}(y)}{\Sigma_{j}(y)}\overrightarrow{\beta}
+\overrightarrow{\gamma}.
\end{equation}

It follows from inequalities (\ref{eq4.8}) and (\ref{eq3.3}) that
\begin{equation}\label{eq4.12}
|\overrightarrow{\gamma}|\leqslant
\sum_{l=1}^{j-1}\dfrac{C|yz|}{R\delta\Sigma_j(y)}\cdot C\nu
\leqslant \dfrac{C|yz|\nu}{R\delta\Sigma_j(y)}.
\end{equation}

In order to estimate $|\overrightarrow{\beta}|$, we introduce two
points $\overline{z}'_{j-1}$ and $z'_{j-1}$ such that
$$\overline{z}'_{j-1}=g_{j-1}^{-1}\left(\sum\limits_{l=1}^{j-1}
\dfrac{\phi_l(y)}{\Sigma_{j-1}(y)}g_{j-1}(h_l(z))\right)$$ and
\begin{equation}\label{eq15}
\overrightarrow{g_{j-1}(y_{j-1})g_{j-1}(z'_{j-1})}=
\overrightarrow{g_{j-1}(\overline{y}_{j-1})g_{j-1}(\overline{z}_{j-1})}.
\end{equation}
Now we put
$$\begin{aligned} &\overrightarrow{\beta}^1=\overrightarrow{\widetilde y_{j-1}\widetilde
z_{j-1}}-\overrightarrow{\widetilde y_{j-1}\widetilde{z'}_{j-1}},\\
&\overrightarrow{\beta}^2=\overrightarrow{\widetilde
y_{j-1}\widetilde{z'}_{j-1}}-\overrightarrow{\widetilde{\overline{y}}_{j-1}\widetilde
{\overline{z}}_{j-1}},\\
&\overrightarrow{\beta}^3=\overrightarrow{\widetilde{\overline{y}}_{j-1}\widetilde
{\overline{z}}_{j-1}}-
\overrightarrow{\widetilde{\overline{y}}_{j-1}\widetilde{\overline{z}'}_{j-1}},\\
&\overrightarrow{\beta}^4=
\overrightarrow{\widetilde{\overline{y}}_{j-1}\widetilde{\overline{z}'}_{j-1}}-
\sum_{l=1}^{j-1}\dfrac{\phi_l(y)}{\Sigma_{j-1}(y)}
\overrightarrow{\widetilde{h_l(y)}\widetilde{h_l(z)}}.
\end{aligned}
$$
Obviously
$\overrightarrow{\beta}=\overrightarrow{\beta}^1+\overrightarrow{\beta}^2
+\overrightarrow{\beta}^3+\overrightarrow{\beta}^4$.

Firstly, $$ \begin{aligned} |\overrightarrow{\beta}^1|=&
\left|\overrightarrow{\widetilde{z'}_{j-1}\widetilde
z_{j-1}}\right|=
(1+\varkappa(\delta))|z'_{j-1}z_{j-1}|\\
=&(1+\varkappa(\delta))\left|\overrightarrow{g_{j-1}(z'_{j-1})g_{j-1}(z_{j-1})}\right|\\
=&(1+\varkappa(\delta))\left|\overrightarrow{g_{j-1}(y_{j-1})g_{j-1}(z_{j-1})}
   -\overrightarrow{g_{j-1}(y_{j-1})g_{j-1}(z'_{j-1})}\right|\\
\text{(by
(\ref{eq15}))}=&(1+\varkappa(\delta))\left|\overrightarrow{g_{j-1}(y_{j-1})g_{j-1}(z_{j-1})}-
\overrightarrow{g_{j-1}(\overline{y}_{j-1})g_{j-1}(\overline{z}_{j-1})}\right|\\
= & (1+\varkappa(\delta))|\overrightarrow{\alpha}_{j-1}|.
\end{aligned}
$$

Secondly, $$\begin{aligned}
|\overrightarrow{\beta}^3|=&\left|\overrightarrow{
\widetilde{\overline{z}'}_{j-1}\widetilde{\overline{z}}_{j-1}}
\right|=
(1+\varkappa(\delta))|\overline{z}'_{j-1}\overline{z}_{j-1}|\\
=& (1+\varkappa(\delta))|\overrightarrow{g_{j-1}(\overline{z}'_{j-1})g_{j-1}(\overline{z}_{j-1})}|\\
 = & (1+\varkappa(\delta))\left|\sum\limits_{l=1}^{j-1}\left(
\dfrac{\phi_l(y)}{\Sigma_{j-1}(y)}-\dfrac{\phi_l(z)}{\Sigma_{j-1}(z)}\right)g_{j-1}(h_l(z))\right|\\
= & (1+\varkappa(\delta))\left|\sum\limits_{l=1}^{j-1}\left(
\dfrac{\phi_l(y)}{\Sigma_{j-1}(y)}-\dfrac{\phi_l(z)}{\Sigma_{j-1}(z)}\right)
\overrightarrow{g_{j-1}(h_1(z))g_{j-1}(h_l(z))}\right| \\
\leqslant &  \dfrac{C|yz|\nu}{\delta R\Sigma_{j-1}(y)}\ \ \
\text{(similar to getting (\ref{eq4.12}))}.
\end{aligned}
$$

Thirdly, we estimate $|\overrightarrow{\beta}^4|$. According to
Lemma 2.4, it follows from (\ref{eq4.?}) and (\ref{eq4.6}) that for
any $1\leqslant l, l_1, l_2\leqslant N_2$
\begin{equation}\label{eq4.*}
\overrightarrow{g_l(h_{l_1}(y))g_l(h_{l_1}(z))} \text{ is
$\varkappa(\delta)$-almost equal to }
\overrightarrow{g_l(h_{l_2}(y))g_l(h_{l_2}(z))},
\end{equation}
and thus
$$\overrightarrow{g_{j-1}(\overline{y}_{j-1})g_{j-1}(\overline{z}'_{j-1})}\text{
is $\varkappa(\delta)$-almost equal to
}\overrightarrow{g_{j-1}(h_{l}(y))g_{j-1}(h_{l}(z))}.$$ Then
according to Corollary \ref{coro2.10}\footnote{When applying
Corollary \ref{coro2.10}, we can assume that $(s_i^{j-1},t_i^{j-1})$
is also an $R$-long $(n,2\delta)$-strainer at $h_j(x_j)$ (see the
beginning of the proof of (\ref{eq4.6}) in Appendix).} and Lemma
2.4,
\begin{equation}\label{eq4.???}\overrightarrow{\widetilde{\overline{y}}_{j-1}\widetilde{\overline{z}'}_{j-1}}
\text{ is $\varkappa(\delta)$-almost equal to }
\overrightarrow{\widetilde{h_{l}(y)}\widetilde{h_{l}(z)}}.
\end{equation} On the other hand, by (\ref{eq4.*})
$$\sum_{l=1}^{j-1}\dfrac{\phi_l(y)}{\Sigma_{j-1}(y)}
\overrightarrow{\widetilde{h_l(y)}\widetilde{h_l(z)}} \text{ is
$\varkappa(\delta)$-almost equal to }
\overrightarrow{\widetilde{h_{l}(y)}\widetilde{h_{l}(z)}}.$$
Therefore it follows that
$$|\overrightarrow{\beta}^4|<\varkappa(\delta)|\widetilde{h_{l}(y)}\widetilde{h_{l}(z)}|=\varkappa(\delta)|yz|.$$

Finally, we estimate $|\overrightarrow{\beta}^2|$. Note that it
follows from (\ref{eq4.???}) that
$|\overrightarrow{\widetilde{\overline{y}}_{j-1}\widetilde{\overline{z}'}_{j-1}}|
<\varkappa(\delta)|yz|$, and thus
$$|\overrightarrow{\widetilde{\overline{y}}_{j-1}\widetilde
{\overline{z}}_{j-1}}|\leqslant |\overrightarrow{\beta}^3|+
|\overrightarrow{\widetilde{\overline{y}}_{j-1}\widetilde{\overline{z}'}_{j-1}}|
<\dfrac{C|yz|\nu}{\delta R\Sigma_{j-1}(y)}+\varkappa(\delta)|yz|.$$
On the other hand, according to Corollary \ref{coro2.10} and Lemma
2.4 it follows from (\ref{eq15}) that
$$
\overrightarrow{\widetilde{y}_{j-1}\widetilde{z'}_{j-1}} \text{ is
$\varkappa(\delta)$-almost equal to }
\overrightarrow{\widetilde{\overline{y}}_{j-1}\widetilde{\overline{z}}_{j-1}}.
$$
Therefore we have
$$ |\overrightarrow{\beta}^2| \leqslant
\varkappa(\delta)|\overrightarrow{\widetilde{\overline{y}}_{j-1}\widetilde{\overline{z}}_{j-1}}|
\leqslant \varkappa(\delta)\left(\dfrac{C|yz|\nu}{\delta
R\Sigma_{j-1}(y)}+\varkappa(\delta)|yz|\right).
$$

Now we can conclude that
$$|\overrightarrow{\beta}|\leqslant|\overrightarrow{\beta}^1|+
|\overrightarrow{\beta}^2|+|\overrightarrow{\beta}^3|+|\overrightarrow{\beta}^4|<
(1+\varkappa(\delta))|\overrightarrow{\alpha}_{j-1}|+\dfrac{C|yz|\nu}{\delta
R\Sigma_{j-1}(y)}+\varkappa(\delta)|yz|.$$ And plugging the
estimates of $|\overrightarrow{\beta}|$ and
$|\overrightarrow{\gamma}|$ (see (\ref{eq4.12})) into
(\ref{eq4.11}), we obtain the {\bf Subclaim} (and thus {\bf the
whole proof is completed}). \hfill$\Box$


\setcounter{equation}{0}

\section{Appendix}

In Appendix, we give the proofs of (2.1.1), Lemma \ref{lemma2.7} and
(\ref{eq4.6}). In the proof of (2.1.1), we will use a result
contained in Lemma 5.6 in \cite{BGP}.
\begin{lemma}\label{lemma5.6}
Let $p,q, r,s\in M$. For sufficiently small $\delta$, if
$|qs|<\delta\cdot\min\{|pq|,|rq|\}$ and $\tilde \angle
pqr>\pi-\delta$, then $|\tilde \angle pqs-\angle
pqs\footnote{$\angle pqs$ is the angle between geodesics $qp$ and
$qs$ at $q$, which is well defined by
$\lim\limits_{x,y\longrightarrow q}\tilde\angle xqy$ with $x\in qp$
and $y\in qs$.}|<\varkappa(\delta)$ and $|\tilde \angle rqs-\angle
rqs|<\varkappa(\delta)$.
\end{lemma}

\noindent{\bf Proof of (2.1.1)}:

According to Lemma \ref{lemma5.6}, (2.1.1) is equivalent to
\begin{equation}\label{eq5.1}
|\angle a_iq_jr_j-\angle a_iq_{j'}r_{j'}|<\varkappa(\delta) \iff
|\angle s_iq_jr_j-\angle s_iq_{j'}r_{j'}|<\varkappa(\delta) \text{
for } i=1,\cdots,n.
\end{equation}

Using the law of cosine, it is not difficult to conclude $$|\tilde
\angle uq_j v-\tilde \angle uq_{j'}v|<\varkappa(\delta) \text{ for }
u\in \{s_i,t_i\}_{i=1}^n \text{ and } v\in \{a_i,b_i\}_{i=1}^n.$$ By
Lemma \ref{lemma5.6} again,
\begin{equation}\label{eq5.2}|\angle
uq_j v- \angle uq_{j'}v|<\varkappa(\delta).\end{equation}

Now we consider spaces of directions at $q_j$, $\Sigma_{q_j}$, with
angle metric. In the situation here, Theorem 9.5 in \cite{BGP}
ensures that $\Sigma_{q_j}$ is $\varkappa(\delta)$-almost isometric
to an $(n-1)$-dimensional unit sphere. Denote by $\bar
a_i\in\Sigma_{q_j}$ (resp. $\bar s_i$ and $\bar r_j$) the directions
of geodesics $q_ja_i$ (resp. $q_js_i$ and $q_jr_j$) for
$i=1,\cdots,n$. Note that
$$|\bar a_i\bar a_{i'}|=\frac\pi2\pm\varkappa(\delta) \text{ and }
|\bar s_i\bar s_{i'}|=\frac\pi2\pm\varkappa(\delta) \text{ for }
i\ne i'.$$ Then it is not difficult to see that inequality
(\ref{eq5.2}) implies (\ref{eq5.1}). \hfill$\Box$

\vskip2mm

\noindent{\bf Proof of Lemma \ref{lemma2.7}}:

We only give the proof for $k=0$ (proofs for other cases are
similar). We first note that
$$
\begin{aligned}  &  |\tilde \angle a_ix_1y_1- \tilde\angle
a_ix_2y_{2}|<\varkappa(\delta) \\
\iff &   |\cos \tilde \angle a_ix_1y_1-\cos \tilde\angle
a_ix_2y_{2}|< \varkappa(\delta)\\
\iff &
\left|\dfrac{|a_ix_1|^2+|x_1y_1|^2-|a_iy_1|^2}{2|a_ix_1|\cdot|x_1y_1|}
 -\dfrac{|a_ix_2|^2+|x_2y_2|^2-|a_iy_2|^2}{2|a_ix_2|\cdot|x_2y_2|}\right|<
\varkappa(\delta) \\
\hskip1cm\iff & \left|\dfrac{|a_ix_1|-|a_iy_1|}{|x_1y_1|}
 -\dfrac{|a_ix_2|-|a_iy_2|}{|x_2y_2|}\right|<
\varkappa(\delta)  \hskip5cm \text{(5.3)}\\
\iff & \left|\dfrac{|a_ix_1|-|a_iy_1|}{|f(x_1)f(y_1)|}
 -\dfrac{|a_ix_2|-|a_iy_2|}{|f(x_2)f(y_2)|}\right|<
\varkappa(\delta) \ \text{($f$ is a $\varkappa(\delta)$-almost
isometry)}.
\end{aligned}
$$
Recall that $f(x)=(|a_1 x|, |a_2 x|, \cdots, |a_n x|)$. Hence
$|\tilde \angle a_ix_1y_1- \tilde\angle a_ix_2y_{2}|<
\varkappa(\delta) $ for $i=1,2,\cdots,n \iff \angle
(\overrightarrow{f(x_1)f(y_1)},\overrightarrow{f(x_2)f(y_2)})<\varkappa(\delta)$.
\hfill$\Box$

\vskip2mm

\noindent{\bf Proof of (\ref{eq4.6})}:

We only give the proof for $k=0$.

We first give an observation that $\{s_i^{j},t_i^{j}\}_{i=1}^n$ is
an $R$-long $(n,C\delta)$-strainer at any $x_l$ for $l=1,\cdots,N_2$
(note that $|x_jx_l|\leqslant N_2R\delta\leqslant NR\delta$ with $N$
depending only on $n$). Without loss of generality, we can assume
that $\{s_i^{j},t_i^{j}\}_{i=1}^n$ is an $R$-long
$(n,\delta)$-strainer at $x_l$, and thus
$\{h(s_i^{j}),h(t_i^{j})\}_{i=1}^n$ is an $R$-long
$(n,2\delta)$-strainer at $h(x_l)$.

Next we note that inequality (\ref{eq4.6}) is equivalent to for any
$1\leqslant j, l_1, l_2\leqslant N_2$
$$|\tilde \angle h(a_i^j)h_{l_1}(y)h_{l_1}(z)- \tilde\angle
h(a_i^j)h_{l_2}(y)h_{l_2}(z)|< \varkappa(\delta).$$ On the other
hand, for $i=1,\cdots,n$ and any $u\in\{s_i^{j},t_i^{j}\}_{i=1}^n$
$$
\begin{aligned} & |\tilde \angle h(a_i^j)h_{l_1}(y)h_{l_1}(z)- \tilde\angle
h(a_i^j)h_{l_2}(y)h_{l_2}(z)|< \varkappa(\delta)\\
\text{(by (2.1.1))}\iff & |\tilde \angle h(u)h_{l_1}(y)h_{l_1}(z)-
\tilde\angle
h(u)h_{l_2}(y)h_{l_2}(z)|< \varkappa(\delta)\ \\
\text{ (obviously)}\Longleftarrow\ & |\tilde \angle
h(u)h_l(y)h_l(z)-\tilde \angle
uyz|< \varkappa(\delta) \text{ for } l=1,\cdots,N_2\\
\text{(by (\ref{lemma5.6}))}\iff & |\angle h(u)h_l(y)h_l(z)-\angle
uyz|<
\varkappa(\delta) \\
(?)\iff & |\angle h(a_i^l)h_{l}(y)h_{l}(z)- \angle
a_i^lyz|< \varkappa(\delta) \hskip4cm\text{(5.4)}\\
\text{(by Lemma \ref{lemma5.6})}\iff& |\tilde \angle
h(a_i^l)h_{l}(y)h_{l}(z)- \tilde\angle
a_i^lyz|< \varkappa(\delta) \\
\text{ (see (5.3))}\iff &
\left|\dfrac{|h(a^l_i)h_l(y)|-|h(a^l_i)h_l(z))|}{|h_l(y)h_l(z)|} -
\dfrac{|a^l_i y| -|a^l_i z| }{ |yz|}\right|<\varkappa(\delta) ,
\end{aligned}
$$
where the last inequality holds because $|h(a^l_i)h_l(y)|=|a^l_iy|$
and $|h(a^l_i)h_l(z)|=|a^l_iz|$ (recall that $h_l=g_l^{-1}\circ
f_l$), and $h_l$ is a $\varkappa(\delta)$-almost isometry.

Hence we only need to verify the third `$\iff$' in (5.4). Similar to
getting inequality (\ref{eq5.2}), we can obtain for any
$v\in\{a_i^l, b_i^l\}_{i=1}^n$
$$|\angle h(u)h_l(y)h(v)-\angle uyv|<\varkappa(\delta).$$
Therefore we can use the same argument as the end of the proof of
(2.1.1) to conclude the third `$\iff$' in (5.4) holds (taking into
account that both $\Sigma_{h_l(y)}$ and $\Sigma_{y}$ are
$\varkappa(\delta)$-almost isometric to $\mathbb{S}^{n-1}$).
\hfill$\Box$

\vskip5mm

{\bf Acknowledgement}: The authors are indebted to Professor Yuri
Burago for his very precious suggestions on the notations,
formulations and structure (and so on) of the present paper.

\end{document}